\documentstyle[12pt]{amsart}
\def\frk{\frak}               

\def\mm{{\frk m}}

\def\opn#1#2{\def#1{\operatorname{#2}}} 
\opn\chara{char}
\opn\length{\ell}
\opn\pd{pd}
\opn\rk{rk}
\opn\projdim{proj\,dim}
\opn\rank{rank}
\opn\depth{depth}
\opn\grade{grade}
\opn\height{height}
\opn\embdim{emb\,dim}
\opn\codim{codim}

\opn\Tr{Tr}
\opn\bigrank{big\,rank}
\opn\superheight{superheight}\opn\lcm{lcm}
\opn\trdeg{tr\,deg}%
\opn\reg{reg}
\opn\lreg{lreg}
\opn\div{div}
\opn\Div{Div}
\opn\cl{cl}
\opn\Cl{Cl}
\opn\Spec{Spec}
\opn\Supp{Supp}
\opn\supp{supp}
\opn\Sing{Sing}
\opn\Ass{Ass}
\opn\Ann{Ann}
\opn\Rad{Rad}
\opn\Soc{Soc}
\opn\Ker{Ker}
\opn\Coker{Coker}
\opn\Im{Im}
\opn\Hom{Hom}
\opn\Tor{Tor}
\opn\Ext{Ext}
\opn\End{End}
\opn\Aut{Aut}
\opn\id{id}

\opn\nat{nat}
\opn\pff{pf}
\opn\Pf{Pf}
\opn\GL{GL}
\opn\SL{SL}
\opn\mod{mod}
\opn\ord{ord}
\opn\aff{aff}
\opn\con{conv}
\opn\relint{relint}
\opn\st{st}
\opn\lk{lk}
\opn\cn{cn}
\opn\core{core}
\opn\vol{vol}
\opn\link{link}
\opn\star{star}
\opn\gr{gr}

\def\pot#1#2{#1[\kern-0.28ex[#2]\kern-0.28ex]}

\opn\dirlim{\underrightarrow{\lim}}
\opn\inivlim{\underleftarrow{\lim}}
\let\union=\cup
\let\sect=\cap
\let\dirsum=\oplus
\let\tensor=\otimes
\let\iso=\cong
\let\Union=\bigcup

\let\Dirsum=\bigoplus

\let\mcone= * 
\let\to=\rightarrow
\let\To=\longrightarrow
\def\Implies{\ifmmode\Longrightarrow \else
     \unskip${}\Longrightarrow{}$\ignorespaces\fi}
\def\implies{\ifmmode\Rightarrow \else
     \unskip${}\Rightarrow{}$\ignorespaces\fi}
\def\iff{\ifmmode\Longleftrightarrow \else
     \unskip${}\Longleftrightarrow{}$\ignorespaces\fi}

\let\:=\colon
\newtheorem{Theorem}{Theorem}[section]
\newtheorem{Lemma}[Theorem]{Lemma}
\newtheorem{Corollary}[Theorem]{Corollary}

\newtheorem{Example}[Theorem]{Example}

\newtheorem{Definition}[Theorem]{Definition}

\let\epsilon\varepsilon
\let\phi=\varphi
\let\kappa=\varkappa
\def\qed{\ifhmode\textqed\fi
   \ifmmode\ifinner\quad\qedsymbol\else\dispqed\fi\fi}
\def\textqed{\unskip\nobreak\penalty50
    \hskip2em\hbox{}\nobreak\hfil\qedsymbol
    \parfillskip=0pt \finalhyphendemerits=0}
\def\dispqed{\rlap{\qquad\qedsymbol}}

\opn\inii{in}
\opn\inim{inm}
\opn\set{set}

\begin{document}

\title{Resolutions by mapping cones}
\author{J\"urgen Herzog and Yukihide Takayama}
\address{J\"urgen Herzog, Fachbereich Mathematik und Informatik,  
Universit\"at-GHS Essen, 45117 Essen, Germany}
\email{juergen.herzog@@uni-essen.de}
\address{Yukihide Takayama, Department of Mathematical Sciences, Ritsumeikan University, 
1-1-1 Nojihigashi, Kusatsu, Shiga 525-8577, Japan}
\email{takayama@@se.ritsumei.ac.jp}

\maketitle

\section*{Introduction}
Many well-known free resolutions arise as  iterated mapping cones. Prominent examples
 are the Eliahou-Kervaire resolution of stable monomial ideals (as noted by Evans and 
Charalambous \cite{EC}), and the Taylor resolution. The idea of the iterated mapping cone
 construction is the following: Let $I\subset R$ be an ideal generated by $f_1,\ldots, f_n$,
  and set $I_j=(f_1,\ldots, f_j)$. Then for $j=1,\ldots, n$ there are exact sequences
\[
0\To R/(I_{j-1}:f_j)\To R/I_{j-1}\To R/I_j\To 0.
\]
Assuming that the free $R$-resolution $F$ of $R/I_{j-1}$ is already known, and a free 
$R$-resolution $G$ of $R/(I_{j-1}:f_j)$ is also known,  one obtains the resolution of  
$R/I_j$ as a mapping cone of a complex homomorphism $\psi\: G\to F$ 
which is a lifting of  $R/(I_{j-1}:f_j)\to R/I_{j-1}$.  Of course one cannot expect 
that such a resolution will be minimal in general. However this construction yields 
an inductive
procedure to compute a resolution  of $R/I$ provided for each $j$, a resolution of 
$R/(I_{j-1}:f_j)$ is known as well as the comparison map $\psi$. 

So it is natural to consider classes of ideals for which the colon ideals $I_{j-1}:f_j$ 
are generated by regular sequences. But even in this nice case it still 
hard to construct the comparison maps. In the first section of this 
paper we therefore restrict ourselves to the case that $I$ is a monomial ideal in 
a polynomial ring, and that the colon ideals in question are generated by subsets of the 
variables. In this case we say that $I$ has linear quotients. At a first glance these 
hypotheses seem to be very restrictive. On the other hand, there are many interesting 
examples of such ideals. All stable and squarefree stable ideals belong to this class, 
as well as all matroidal ideals (see Section \ref{section:linearq} for the definitions).

It is easy to see that $I$ has linear quotients, if and only if the first syzygy module 
of $I$ has a quadratic Gr\"obner basis, which, in case of a Stanley  ideal $I_\Delta$ 
attached to the              simplicial complex $\Delta$, is equivalent to saying that 
the Alexander dual  $\Delta^*$ of $\Delta$ is nonpure shellable. This fact was communicated
 to us by Sk\"oldberg.  It is also easy to see that $I$ has linear quotients if and only 
 if $I$ satisfies condition $(4.1)$ of Batzies and Welker \cite{BW}, a condition  which the 
 authors call shellable.  

It is clear that our approach   only requires to describe the comparison maps in order to 
compute explicit free resolutions  of ideals with linear quotients as iterated mapping cones.  
Our description  of the comparison maps is modeled 
after Eliahou and Kervaire and is based on decomposition functions. A function $g$ which 
assigns to each monomial in $I$ (in a natural way) a monomial generator of $I$, 
see \ref{lemma:decomposition}, is called a decomposition function. If it satisfies 
a certain additional condition which is described in 
Definition \ref{def:regular}, then we call it regular. Stable and squarefree stable ideals 
have regular 
decomposition functions, but also matroidal ideals as we show in 
Theorem \ref{theorem:matroid}. The main result of Section \ref{section:linearq} 
however is Theorem \ref{theorem:main} in which we give an explicit resolution  of 
all ideals with linear quotients which admit a regular decomposition function. 
These include then of course also matroidal ideals for which explicite resolutions 
in different terms are already known by Reiner and Welker \cite{RW},  
and Novik, Postnikov and Sturmfels \cite{NPS}.

In the second part of this paper we ask ourselves under which circumstances a mapping cone can 
be given the structure of a DG algebra. The natural way of doing 
this, is to assume that $F$     
(notation as above) is already a DG algebra, that $G$ is a DG $F$-module, 
and that the mapping cone of $\psi$ is a trivial extension of $F$ by $G$ 
in the category of DG modules. This idea was first used by Levin and Avramov  \cite{LA} 
in order to compute the Poincar\'e series of Gorenstein algebra modulo its socle. 
In Section \ref{section:trivialex} we analyze under which conditions the mapping 
cone can be given the structure of a trivial extension and in the final 
Section \ref{section:koszultype} we describe cases for which resolutions 
which are constructed as iterated mapping cones admit a DG algebra resolution.  
We call such resolutions of Koszul type and show in Example~\ref{example:taylor} that the
known DG algebra structure (see \cite{F} and \cite{G}) on the Taylor complex is 
of Koszul type, as well as the resolution of an almost complete intersection which 
is directly linked to a complete intersection, see Example~\ref{example:almost}.

\section{Monomial Ideals with linear quotients
				  \label{section:linearq}}
Let $K$ be a field, $R=K[x_1,\ldots,x_n]$ be the polynomial ring in $n$ indeterminates, 
and $I\subset R$ a monomial ideal. The unique minimal set of monomial generators of $I$ 
will be denoted by $G(I)$.  The ideal $I$ is said to have {\em linear quotients} if for 
some order $u_1,\ldots,u_m$ of the elements of  $G(I)$ and all $j=1,\ldots,m$  the colon 
ideals $(u_1,\ldots,u_{j-1}):u_j$ are generated by a subset of  $\{x_1,\ldots,x_n\}$. 

We define 
\[
\set(u_j)=\{k\in[n]\: x_k\in (u_1,\ldots,u_{j-1}):u_j\}\quad \text{for} 
\quad j=1,\ldots, m.
\]
 
\begin{Example}
\label{example:strongly}
{\em According to Eliahou-Kervaire \cite{EK}, a monomial ideal $I$ is called 
{\em stable} if for all $u\in G(I)$ and all $i\leq m(u)$ one has that 
$x_i(u/x_{m(u))})\in I$. Here 
$m(u)=\max\{i\: i\in\supp(u)\}$, and $\supp(u)=\{i\: x_i\quad\text{divides}\quad u\}$.

Let $G(I)=\{u_1,\ldots,u_m\}$, where $u_1>u_2>\cdots > u_m$ in the reverse degree 
lexicographical order with regard to $x_1 > x_2 > \ldots > x_n$. It is easy to see 
that $I$ has linear quotient for this order  of the
generators, and that $\set(u)=\{1,\ldots,m(u)-1\}$ for all $u\in G(I)$.}
\end{Example}

\begin{Example}
\label{example:sqstable}
{\em In \cite{AHH} a squarefree monomial ideal $I$ is called 
{\em squarefree stable} if for all $u\in G(I)$ and all $i<m(u)$ with 
$i\not \in\supp(u)$ one has that $x_i(u/x_{m(u))})\in I$. 
With respect to the reverse degree lexicographical order on the generators, 
$I$ has linear quotients and $\set(u)=\{i\: i<m(u), i\not\in\supp(u)\}$. } 
\end{Example}

We  now want to analyze more carefully when a squarefree monomial ideal $I$ has 
linear quotients. Let  $\Delta$ the corresponding simplicial complex 
on the vertex set $[n]=\{1,\ldots,n\}$, so that $I=I_{\Delta}$. 
The {\em Alexander dual of $\Delta$}  is the simplicial complex 
$\Delta^*=\{\sigma\in [n]\: [n]\setminus \sigma\not\in \Delta\}$. 
The following two statements are almost tautologically equivalent: 
(i) $I_\Delta$ has linear quotients, 
(ii) $\Delta^*$ is shellable in the non-pure sense of Bj\"orner and Wachs \cite{BjW}.

In particular, it follows that if the simplicial complex 
$\Gamma$ generated by 
\[{\cal B} = \{\supp(u) \: u\in G(I_{\Delta})\}
\] 
is a matroid (in which case we say that $I_{\Delta}$ is {\em matroidal}), then
$I_{\Delta}$ has linear quotients. 
In fact, if $\Gamma$ is  a matroid, then ${\cal B}$ is the matroid basis, and 
$\Delta^*$ is the dual matroid of $\Gamma$ with basis  
$\{[n]\backslash\supp(u) \: u\in G(I_{\Delta})\}$. 
On the other hand, it is known that all matroids are shellable, see \cite{Bj}.

For the convenience of the reader we show directly that,
for a simplicial complex $\Delta$,
the Stanley-Reisner ideal $I_{\Delta}$
has linear quotients, if it is  matroidal.
If $I_{\Delta}$ is matroidal, then
for all $u,v\in G(I_\Delta)$ and all $i\in\supp(u)\setminus \supp(v)$ 
there exists $j\in \supp(v)\setminus \supp(u)$ such that 
$x_j(u/x_i)\in I_{\Delta}$. 
We now claim that the generators of $I_\Delta$ in reverse lexicographical order have 
linear quotients.

We will prove  a slightly more general result  from which this claim will follow. 
Let $u=x_1^{a_1}\cdots x_n^{a_n}$ be a monomial. We set $\nu_i(u)=a_i$ for $i=1,\ldots,n$.

\begin{Lemma}
\label{moregeneral}
Let $I$ be a monomial ideal for which all generators have the same degree. 
Suppose $I$ satisfies the following exchange property: 
\begin{itemize}
\item[]For all $u,v\in G(I)$ and all $i$ with $\nu_i(u)>\nu_i(v)$, there exists 
an integer $j$ with $\nu_j(v)>\nu_j(u)$ such that $x_j(u/x_i)\in G(I)$. 
\end{itemize}
Then $I$ has linear quotients with respect to the reverse lexicographic order 
of the generators.
\end{Lemma}

Note that a squarefree monomial ideal is matroidal if and only if it satisfies 
the exchange property of Lemma \ref{moregeneral}. Blum \cite{Bl} has shown that 
ideals whose generators correspond to the basis of a polymatroid satisfy the exchange property 
of Lemma~\ref{moregeneral}.

\begin{pf}[Proof of \ref{moregeneral}]
Let $u\in G(I)$, and let $J$ be the ideal generated 
by all $v\in G(I)$ with $v>u$ (in the reverse lexicographical order). Then 
\[
J:u=(v/[v,u]\: v\in J ),
\]
where $[v,u]$ denotes the greatest common divisor of $v$ and $u$. 
Thus in order to prove that $J:u$ is generated by monomials  of degree $1$, 
we have to show that for each $v>u$ there exists $x_j\in J:u$ such that $x_j$ divides 
$v/[v,u]$.  

In fact, let
$u=x_1^{a_1}\cdots x_n^{a_n}$ and 
$v=x_1^{b_1}\cdots x_n^{b_n}$. Since $v>u$, 
there exists an integer $i$ such that $a_k=b_k$ for $k=i+1,\ldots,n$,  
and $a_i > b_i$, and hence an integer $j$ with $b_j>a_j$ such that  
$u'=x_j(u/x_{i})\in I$. 
Since $j < i$, we see that $u'\in J$, and from the equation 
$x_iu'=x_ju$ we deduce that $x_j\in J:u$. Finally, 
since $\nu_j(v/[u,v])=b_j-\min\{b_j,a_j\}=b_j-a_j>0$, 
we have that $x_j$ divides $v/[v,u]$.
\end{pf}
\medskip

Since for a matroid 
$\Delta^*$, 
the Stanley-Reisner ideal $I_\Delta$ has linear quotients with respect to 
the reverse lexicographical order of the generators, it follows that 
\[
\set(u)\subset \{i\: i<m(u), i\not\in\supp(u)\}
\] 
for all $u\in G(I_\Delta)$. More precisely, we have 
\[
\set(u)=\{i \ : \supp(v)\setminus \supp(u)=\{i\}
          \quad\text{for some}\quad v\in G(I_\Delta), v>u\}.
\]  

\begin{Example}
\label{example:transversal}
{\em Let $A_1,\ldots,A_r$ be non-empty subsets of $[n]$.  The collection of subsets
$\{i_1,\ldots,i_r\}$ of $[n]$ with $i_j\in A_j$ for $j=1,\ldots,r$ 
and $i_j\neq i_k$ for $j\neq k$, is the basis of a matriod, called {\em transversal}. 
Let $I$ be the squarefree monomial ideal whose generators correspond to the basis of 
this transversal matriod. Let $u\in G(I)$, $u=x_{i_1}\cdots x_{i_r}$. 
The above description of $\set(u)$ yields in this case
\[
\set(u)=\Union_{j=1}^r\{k\in A_j\: k\leq j\}\setminus \supp(u).
\]
}
\end{Example}
\medskip

Let $\psi\: A \to B$ be a complex homomorphism. 
Recall that the {\em mapping cone of $\psi$} is the complex $C(\psi)$ with 
$C(\psi)_i=B_i\dirsum A_{i-1}$ for all $i$, and chain map 
$d$ with 
$d_i\: C(\psi)_i\to C(\psi)_{i-1}$,
$d_i(b,a)=(\psi(a)+\partial(b), -\partial(a))$.  

We want to apply this concept in the following situation: 
Suppose that $I$ has linear quotients with respect the order $u_1,\ldots, u_m$ of 
the generators of $I$. Set $I_j=(u_1,\ldots,u_j)$ and $L_j=(u_1,\ldots,u_j):u_{j+1}$, 
then, since $I_{j+1}/I_j\iso R/L_j$, we get the exact sequences
\[
0\To R/L_{j}\To R/I_j\To R/I_{j+1}\To 0.
\]
The homomorphism $R/L_{j+1}\to R/I_j$ is multiplication with $u_{j+1}$. Let $F^{(j)}$ 
be a graded  free resolution of $R/I_j$, $K^{(j)}$ the Koszul complex for the regular
sequence $x_{k_1},\ldots, x_{k_l}$ with $k_i\in \set(u_{j+1})$,  
and $\psi^{(j)}\: K^{(j)}\to F^{(j)}$ a graded complex homomorphism lifting 
$R/L_{j}\to R/I_j$. Then the mapping cone $C(\psi^{(j)})$ of $\psi^{(j)}$ yields 
a free resolution of $R/I_{j+1}$. Thus by iterated mapping cones we obtain step by step 
a graded  free resolution of $R/I$.  

\begin{Lemma}
\label{lemma:basis}
Suppose $\deg u_1\leq \deg u_2\leq \cdots\leq \deg u_m$. Then the iterated mapping cone 
$F$, derived from the sequence $u_1,\ldots,u_m$, 
is a minimal graded  free resolution of $R/I$, and for all $i>0$, the symbols
\[
f(\sigma; u)\quad \text{with} \quad  u\in G(I), 
           \quad \sigma\subset \set(u), \quad |\sigma|=i-1
\]
form a  homogeneous basis of the $R$-module $F_i$. Moreover, $\deg f(\sigma;u)= |\sigma|+\deg 
u$.
\end{Lemma}

We remark that if $I$ has linear quotients with respect to $u_1,\ldots, u_m$, 
then this does not necessarily imply that 
$\deg u_1\leq \deg u_2\leq \cdots\leq \deg u_m$. 
In fact, $I=(x_1x_2, x_2x_3x_4, x_1x_3)$ has linear quotients for the given order 
of the generators.

\begin{pf}[Proof of   Lemma~\ref{lemma:basis}] 
We prove by induction on $j$ that $F^{(j)}$ is a minimal free resolution of $R/I_j$, 
and that $F^{(j)}$ has a basis as asserted. For $j=1$, the assertion is trivial. 
In homological degree $i-1$ the Koszul complex $K^{(j)}$ has the $R$-basis 
$e_{\sigma}=e_{j_1}\wedge\ldots\wedge e_{j_{i-1}}$, 
where $\sigma=\{j_1<j_2<\cdots<j_{i-1}\}\subset\set(u_{j+1})$. 
Since $F^{(j+1)}_i=F^{(j)}_i
\dirsum K^{(j)}_{i-1}$, we obtain the desired basis from the induction hypothesis 
if we identify the elements $e_\sigma$ with $f(\sigma;u_{j+1})$. 

In order to show that $F^{(j+1)}$ is a minimal free resolution, it suffices to show that 
$\Im(\psi^{(j)})\subset \mm F^{(j)}$. Let $f(\sigma;u_{j+1})\in  K^{(j)}_{i-1}$ 
and $\psi^{(j)}(f(\sigma;u_{j+1}))=\sum_{i=1}^{j}\sum_{\tau}a_{\tau, i}f(\tau;u_i)$. 
Since $|\tau|=|\sigma|-1$ and $\deg u_{j+1}\geq \deg u_i$ for all $i=1,\ldots, j$ 
it follows that $\deg  f(\sigma;u_{j+1})>\deg f(\tau;u_i)$   for all $\tau$ and $i$, 
and so $\deg a_{\tau,i}>0$ for all $\tau$ and $i$.\end{pf}

\begin{Corollary}
\label{cor:betti}
The bigraded Poincar\'e series of an ideal with linear quotients is given by
\[
P_{R/I}(s,t)= 1+\sum_{u\in G(I)}(1+s)^{|\set(u)|}st^{\deg u}.
\]
\end{Corollary}

Next we want to describe the chain maps of the graded minimal free resolution of an ideal 
with linear quotients as explicitly as possible. It will turn out that the maps are 
described similarly as in the Eliahou-Kervaire resolution \cite{EK} provided we impose 
some extra condition on the linear quotients.

Let $I$ have linear quotients with respect to the sequence of generators 
$u_1,\ldots,u_m$, and set as before $I_j=(u_1,\ldots,u_j)$ for $j=1,\ldots, m$.  
Let $M(I)$ be the set of all monomials in $I$. 
The  map $g\: M(I)\to G(I)$ is defined as follows: we set $g(u)=u_j$, if $j$ 
is the smallest number such that $u\in I_j$.

\begin{Lemma}
\label{lemma:decomposition}
{\em (a)} For all $u\in M(I)$ one has $u=g(u)c(u)$ with some complementary factor $c(u)$, 
and  $\set(g(u))\sect\supp(c(u))=\emptyset$. 

{\em (b)} Let $u\in M(I)$, $u=vw$ with $v\in G(I)$ and $\set(v)\sect \supp(w)=\emptyset$. 
Then $v=g(u)$.
\end{Lemma}

Notice that any function $M(I)\to G(I)$ satisfying Lemma~\ref{lemma:decomposition}(a) is 
uniquely determined because of Lemma~\ref{lemma:decomposition}(b). 
We call $g$ the {\em decomposition function} of $I$. 

\begin{pf}[Proof of Lemma~\ref{lemma:decomposition}] 
(a) Suppose $g(u)=u_j$. Since $u\in I_j$ it is a multiple of some $u_i$ with $i\leq j$. 
If $i<j$, then  $u\in I_i$, a contradiction. This shows that $g(u)$ divides $u$, i.e., 
$u=g(u)c(u)$ for some $c(u)$. Suppose $\set(g(u))\sect \supp(c(u))\neq \emptyset$, 
and let $i\in \set(g(u))\sect \supp(c(u))$. Then $u=(x_iu_j)(c(u)/x_i)\in I_{j-1}$, 
a contradiction. 

(b) Let $u\in M(I)$. Suppose there exist $u_i, u_j\in G(I)$ such that $u=u_iv_i=u_jv_j$ 
for some monomials $v_i$ and $v_j$,
and $\set(u_i)\sect\supp(v_i)=\emptyset =\set(u_j)\sect\supp(v_j)$. 
We may assume that $i<j$. Then the equation $u_iv_i=u_jv_j$ implies that 
$v_j\in I_{j-1}:u_j$. Hence there exists $k\in \set(u_j)$ such that $x_k|v_j$. 
In other words, we have $k\in \set(u_j)\sect\supp(v_j)$, a contradiction.\end{pf}

The following properties of the decomposition function will be needed later 

\begin{Lemma}
\label{lemma:needed}
Let $u, v\in M(I)$. Then $g(uv)=g(u)$ if and only if 
\[
\set(g(u))\sect \supp(v)=\emptyset.
\]
\end{Lemma}

\begin{pf}
Since $u=g(u)c(u)$, we have $uv=g(u)c(u)v$. Thus if $\set(g(u))\sect\supp(v)=\emptyset$, 
it also follows that $\set(g(u))\sect\supp(c(u)v)=\emptyset$. 
Because of the uniqueness of the decomposition function we conclude that $g(uv)=g(u)$.

Conversely, suppose that $g(uv)=g(u)$. Then $c(u)v=c(uv)$, 
and so $\supp(v)\subset \supp(c(uv))$. Hence, since $\supp(c(uv))$ and $\set(g(uv))\bf $ 
are disjoint sets, $\supp(v)$ and $\set(g(uv))$ are disjoint, too. This yields the 
assertion,
 since $g(u)=g(uv)$.\end{pf}

\begin{Definition}
\label{def:regular}
{\em We say that the decomposition function $g\: M(I)\to G(I)$ is {\em regular}, if
$\set(g(x_su))\subset\set(u)$ for all $s\in\set(u)$ and $u\in G(I)$.

}
\end{Definition} 

Unfortunately the decomposition function for an ideal with linear quotients is not always 
regular. For example, consider 
$I=(x_{2}x_{4}, x_{1}x_{2}, x_{1}x_{3})$.  
Then with respect to the given order of the generators, $I$ has linear quotients. 
One checks that $\set(x_{1}x_{3})=\{2\}$, and that $\set(g(x_2(x_{1}x_{3})))=\{4\}$.

On the other hand it is obvious that stable and squarefree stable ideals have regular 
decomposition functions with respect to the reverse degree lexicographic order. 
Another  large class of squarefree 
ideals with regular decomposition function is given by

\begin{Theorem}
\label{theorem:matroid}
Let $I$ be the Stanley-Reisner ideal of a matroid. Then $I$ has a regular decomposition 
function.
\end{Theorem}

\begin{pf}
Let $G(I) = \{u_1,\ldots, u_n\}$ with  $u_1 > \cdots > u_n$ in the reverse 
lexicographic order. We will set $u = u_n$ for convenience.
Then 
\[
     \set(u) = \{ i : \supp(u_k)\setminus \supp(u) = \{i\}, 
		\mbox{ for some } u_k > u\}.
\]
Take an arbitrary element $i\in \set(u)$, then we have 
\[
     \set(g(x_iu)) = \{ j : \supp(u_l)\setminus \supp(g(x_iu)) = \{j\}, 
		\mbox{ for some } u_l > g(x_iu)\}.
\]
Now we will prove that $g$ is regular, namely that $\set(g(x_i u)) \subset \set(u)$.

We first note the following:
\begin{equation}
\label{one}
     g(x_iu) = \frac{x_iu}{x_{j(i)}}
	\qquad\mbox{for some}\quad j(i)\in \supp(u), \quad i < j(i),
\end{equation}
and for arbitrary $j\in \set(g(x_iu))$, 
\begin{equation}
\label{two}
     \{j\} = \supp(u_l)\setminus ((\supp(u)\cup\{i\})\setminus j(i)).
\end{equation}
Notice that, since $i\in \supp(g(x_iu))$ by (\ref{two}), we have $i\not=j$.
In the following, we will prove that $j \in \set(u)$.

Now we first consider the case of $i\not\in \supp(u_l)$.
From (\ref{two}), we have 
\[
	\{j\} = \supp(u_l) \setminus (\supp(u)\setminus j(i)).
\]
Now assume that $j(i)\in \supp(u_l)$. 
Then we must have $j = j(i)$, and $\supp(u_l) = \supp(u)$, a contradiction.
Thus we have $j(i)\not\in \supp(u_l)$ and 
\[
	\{j\} = \supp(u_l) \setminus \supp(u).
\]
Since $u_l > g(x_iu) > u$, this means that $j\in \set(u)$.

Next we consider the case of $i\in \supp(u_l)$.
From (\ref{two}), we have 
\begin{equation}
\label{three}
	\{i, j\} = \supp(u_l) \setminus (\supp(u)\setminus j(i)).
\end{equation}
Now we have $j(i)\not\in \supp(u_l)$.
In fact, assume that $j(i)\in \supp(u_l)$. Then we have 
$j(i) \in \{i, j\}$ by  (\ref{three}). Moreover, $j(i) = j$ since $i< j(i)$.
But then we have $x_iu = x_{j(i)}g(x_iu) =  x_jg(x_iu)$ by (\ref{one}).
Applying the decomposition function $g$, we obtain
$g(x_iu) = g(x_jg(x_iu))$. This contradicts with the assumption
$j \in \set(g(x_iu))$. Thus we have 
\begin{equation}
\label{four}
	\{i, j\} = \supp(u_l) \setminus \supp(u).
\end{equation}
Since $\{\supp(u_i)\}_i$ is a basis of a matroid, there exists
some $k\in \supp(u)\setminus\supp(u_l)$ such that 
$(x_k/x_i)u_l\in G(I)$. We denote this element by $u_p$.
Moreover we have 
\begin{eqnarray*}
\{j\} & = & \supp(u_l)\setminus (\supp(u)\cup\{i\}\setminus j(i)) \\
      & = & \supp(u_l)\setminus (\supp(u)\cup\{i\})\\
      & = & (\supp(u_l)\setminus i) \setminus \supp(u)\\
      & = & (\supp(u_l)\setminus i)\cup\{k\} \setminus \supp(u)\\
      & = & \supp(u_p)\setminus  \supp(u).
\end{eqnarray*}
If $u_p > u$, this equation means $j\in \set(u)$, and we are done.

Now in the rest of the proof, we assume $u_p < u$. Then, since $u_p =
(x_k/x_i)u_l$ and $u_l > u$, we must have $i< k$.
Assume that $k = j(i)$. Since all bases of a matroid have the 
same cardinality  and $|\supp(u_l)\setminus\supp(u)| = 2$
by (\ref{four}), there must be some element $\zeta\in\supp(u)$ such that 
\[
	\{\zeta, j(i)\} = \supp(u)\setminus\supp(u_l).
\]
Then we have 
\[
	u_l = \frac{x_ix_ju}{x_kx_{\zeta}},
	\qquad
	u_p = \frac{x_ku_l}{x_i} = \frac{x_ju}{x_{\zeta}}
\]
and $u_p < u$ implies $j > \zeta$. Moreover, 
$g(x_iu) = (x_i/x_{j(i)})u = (x_i/x_k)u$, by (\ref{one}). Thus we have 
$u_l = (x_j/x_{\zeta})g(x_iu)$ and 
$u_l < g(x_iu)$. But this contradicts with the 
assumption that $u_l > g(x_iu)$. Consequently, we must have 
$k\not= j(i)$.

Now we have 
\begin{equation}
\label{five}
	\{k, j(i)\} = \supp(u)\setminus\supp(u_l).
\end{equation}
Since 
\[
u_l = \frac{x_jx_iu}{x_kx_{j(i)}} > g(x_iu) = \frac{x_iu}{x_{j(i)}}, 
\]
we have $j < k$. Since moreover we have $i<j(i)$ and $i<k$, there five 
possibilities of total order on
$i$, $j$, $k$ and $j(i)$ left, namely
(i) $i< j < k < j(i)$, (ii) $i < j < j(i) < k$, 
(iii) $j < i < k < j(i)$, (iv) $j < i < j(i) < k$, 
or (v) $i< j(i) < j < k$.
It is easy to check that only in the last case one has $u_p < u$.
So we assume $i< j(i) < j < k$ in the following.

Since $\{\supp(u_i)\}_i$ is the basis of a matroid, 
we have, for $k$, either $(x_i/x_k)u\in G(I)$ or $(x_j/x_k)u\in G(I)$,
by (\ref{four}).
If $(x_j/x_k)u\in G(I)$, then 
\[
	\{j\} = \supp\left(\frac{x_ju}{x_k}\right)
		\setminus\supp(u),
\]
and, since $j < k$, we get $(u_j/x_k)u > u$. Hence we have $j\in \set(u)$.
Now assume that $(x_j/x_k)u\not\in G(I)$. Then we must have 
$(x_i/x_k)u\in G(I)$. Since $i < k$, we get $(x_i/x_k)u > u$.
But $g(x_iu) = (x_i/x_{j(i)})u < (x_i/x_k)u$ since $j(i) < k$,
a contradiction. Hence this case does not happen.\end{pf}

\begin{Lemma}
\label{lemma:exchange}
If $g\: M(I)\to G(I)$ is a regular decomposition function, then 
\[
g(x_sg(x_tu))=g(x_tg(x_su))\quad\text{for all}\quad 
u\in M(I)\quad \text{and all}\quad s,t\in \set(u).
\]
\end{Lemma}

\begin{pf}
If $s\not\in \set(g(x_tu))$, then Lemma~\ref{lemma:needed} implies that 
$g(x_{s}g(x_{t}u))=g(g(x_{t}u))=g(x_{t}u)$, and if $s\in \set(g(x_{t}u))$, then 
$\set(g(x_sg(x_tu)))\subset \set(g(x_tu))$ since $g$ is regular. Thus in any case, 
it follows that  $\set(g(x_sg(x_tu)))\subset \set(g(x_tu))$. Hence 
\[
\set(g(x_sg(x_tu)))\sect \supp(c(x_tu))=\emptyset,
\] 
so that by Lemma~\ref{lemma:needed} again we have $g(x_{s}g(x_{t}u)c(x_{t}u)) =
g(x_{s}g(x_{t}u))$.
Therefore, by the uniqueness of the decomposition function 
(see Lemma~\ref{lemma:decomposition}(b)), 
 the equation $x_sx_tu=(x_sg(x_tu))c(x_tu)$ yields $g(x_sg(x_tu))=g(x_sx_tu)$. 
 This implies the assertion.\end{pf}

The exchange property of the decomposition function in Lemma~\ref{lemma:exchange} is 
weaker than the regularity, as is demonstrated by the following example: The ideal 
$I=(x_1x_3, x_2x_3,\\ x_1x_5, x_3x_4, x_4x_5)$ has linear quotients with respect to 
the given order of the generators. One checks that the exchange property holds. 
But since $\set(x_4x_5)=\{1,3\}$, and  $\set(g(x_3(x_4x_5)))=\{1,2\}$, $g$ is not 
regular.

The following theorem is the main result of this section. It generalizes the theorem of 
Eliahou-Kervaire (\cite{EK}) for stable ideals and that of Aramova-Herzog-Hibi (\cite{AHH}) 
for squarefree stable ideals. For convenience, and to avoid unnecessary distinctions, 
we  extend the definition introduced in Lemma~\ref{lemma:basis} and 
set $f(\sigma;u)=0$ if $\sigma\not\subset \set(u)$.  

\begin{Theorem}
\label{theorem:main}
Let $I$ be a monomial ideal with linear quotients, and $F$ the graded minimal free 
resolution of $R/I$. Suppose the decomposition function $g\: M(I)\to G(I)$ is regular.
Then the chain map $\partial$ of $F$ is given by
\[
\partial(f(\sigma;u))=-\sum_{t\in\sigma}(-1)^{\alpha(\sigma;t)}x_tf(\sigma\backslash t;u)
+\sum_{t\in \sigma}(-1)^{\alpha(\sigma;t)}\frac{x_tu}{g(x_tu)}f(\sigma\backslash t;g(x_tu)),
\]
if $\sigma\neq \emptyset$, and
\[
\partial(f(\emptyset;u))=u\quad \text{otherwise}.
\]
Here $\alpha(\sigma;t)=|\{s\in\sigma\: s<t\}|$.
\end{Theorem}

\begin{pf}
Let $I$ have linear quotients with respect to the sequence $u_1,\ldots, u_m$. 
We show by induction on $j$, that $F^{(j)}$ has the desired chain map. 
For $j=1$, the assertion is trivial. 
Since $F^{(j+1)}$ is the mapping cone of $\psi^{(j)}\: K^{(j)}\to F^{(j)}$ 
it follows that $F^{(j)}$ is a subcomplex of $F^{(j+1)}$ and it suffices to check 
the formula for the chain map on the basis elements $f(\sigma;u_{j+1})$. 
By the definition of the mapping cone of $\psi^{(j)}$ we have 
$\partial(f(\sigma;u_{j+1}))=-\partial_1(f(\sigma;u_{j+1}))+\psi^{(j)}(f(\sigma;u_{j+1}))$,
where $\partial_1$ is the chain map of the Koszul complex $K^{(j)}$. 
Thus in order to prove the asserted formula it remains to show 
that we can define $\psi^{(j)}$ as 
\[
\psi^{(j)}(f(\sigma;u_{j+1}))
= \sum_{t\in \sigma}(-1)^{\alpha(\sigma;t)}\frac{x_tu_{j+1}}{g(x_tu_{j+1})}
f(\sigma\backslash t;g(x_tu_{j+1})),
\]
if $\sigma\neq \emptyset$, and $\psi^{(j)}(f(\emptyset; u_{j+1}))=u_{j+1}$, otherwise. 

To verify this we must prove that  $\psi^{(j)}\circ \partial_1=\partial\circ \psi^{(j)}$. 
In order to simplify notation we set $u=u_{j+1}$ and $\psi=\psi^{(j)}$.
Then for $t\in \set(u)$ we have 
\begin{eqnarray*}
(\psi\circ \partial_1)(f(\{t\};u))&=& \psi(x_tf(\emptyset; u))\\
&=&x_tu,
\end{eqnarray*}
while on the other hand
\begin{eqnarray*}
(\partial\circ \psi)(f(\{t\};u))&=&\partial\left(\frac{x_tu}{g(x_tu)}
f(\emptyset;g(x_tu)\right)\\
=\frac{x_tu}{g(x_tu)}g(x_tu)&=& x_1u.
\end{eqnarray*}
Now let $\sigma\subset \set(u)$ with $|\sigma|\geq 2$. Then 
\begin{eqnarray*}
(\psi\circ \partial_1)(f(\sigma;u))
   &=& \sum_{t\in \sigma}(-1)^{\alpha(\sigma;t)}x_t\psi(f(\sigma\backslash t;u))\\
   &=&\sum_{t\in \sigma}(-1)^{\alpha(\sigma;t)}
            x_t\left(
			    \sum_{s\in\sigma\backslash t}
                (-1)^{\alpha(\sigma\backslash t;s)}\frac{x_su}{g(x_su)}
                f(\sigma\backslash \{s,t\};g(x_su)
				\right)\\
   &=&\sum_{t\in\sigma}\sum_{s\in\sigma\setminus t\atop s<t}
       (-1)^{\alpha(\sigma;t)+\alpha(\sigma;s)}
	   \frac{x_tx_su}{g(x_su)}f(\sigma\backslash\{s,t\}; g(x_su))\\
   &-& \sum_{t\in\sigma}\sum_{s\in\sigma\setminus t\atop s>t}
          (-1)^{\alpha(\sigma;t)+\alpha(\sigma;s)}\frac{x_tx_su}{g(x_su)}
          f(\sigma\backslash\{s,t\}; g(x_su)).
\end{eqnarray*}
Exchanging the role of $s$ and $t$ in the second sum, we obtain
\begin{eqnarray}
\label{a}
(\psi\circ \partial_1)(f(\sigma;u))
&=&\sum_{t\in\sigma}\sum_{s\in\sigma\setminus t\atop s<t}
(-1)^{\alpha(\sigma;t)+\alpha(\sigma;s)}\frac{x_tx_su}{g(x_su)}f(\sigma\backslash\{s,t\}; 
g(x_su))\\
&-&
\sum_{s\in\sigma}\sum_{t\in\sigma\setminus s\atop      
s<t}(-1)^{\alpha(\sigma;t)+\alpha(\sigma;s)}\frac{x_tx_su}{g(x_tu)}
f(\sigma\backslash\{s,t\}; g(x_tu)).\nonumber
\end{eqnarray}
On the other hand we have 
\begin{eqnarray}
\label{b}
(\partial\circ \psi)(f(\sigma;u))&=&
\sum_{t\in\sigma}(-1)^{\alpha(\sigma;t)}\frac{x_tu}{g(x_tu)}
\partial(f(\sigma\backslash t;g(x_tu)),
\end{eqnarray}
and 
\begin{eqnarray*}
\partial(f(\sigma\backslash t;g(x_tu)))=
&-&\sum_{s\in\sigma\setminus t}(-1)^{\alpha(\sigma;s)}x_sf(\sigma\backslash\{s,t\}; g(x_tu))\\
&+& \sum_{s\in\sigma\setminus 
t}(-1)^{\alpha(\sigma\;s)}\frac{x_sg(x_tu)}{g(x_sg(x_tu))}f(\sigma\backslash\{s,t\}; 
g(x_sg(x_tu))
\end{eqnarray*}
Before we continue our calculation we notice that it may happen that 
$\sigma\setminus t\not
\in \set(g(x_tu))$, in which case $f(\sigma\backslash t;g(x_tu))=0$, by convention. 
Thus the right hand side of the equation should also  be zero. 

In fact, let $s\in\sigma\setminus t$. If $\sigma\setminus\{s,t\}\not\subset \set(g(x_tu))$, 
then $\sigma\setminus\{s,t\}\not\subset \set(g(x_sg(x_tu)))$, 
since $g$ is regular, and so the corresponding summands are zero. 
Otherwise, $\sigma\setminus\{s,t\}\subset \set(g(x_tu))$. 
But then $s\not\in \set(g(x_tu))$, so that $g(x_s(g(x_tu)))=g(x_tu)$, 
by Lemma~\ref{lemma:needed}. Therefore,
\[
\frac{x_sg(x_tu)}{g(x_sg(x_tu))}f(\sigma\backslash\{s,t\};g(x_sg(x_tu)))=
x_sf(\sigma\backslash\{s,t\};g(x_tu)).
\]
Hence we see that the summands on the right hand side of the equation are either zero 
or cancel each other, as we wanted to show. 

Now continuing with our calculation we get 
\begin{eqnarray*}
\partial(f(\sigma\backslash t;g(x_tu)))
=&-&\sum_{s\in\sigma\setminus t\atop s<t}(-1)^{\alpha(\sigma;s)}x_sf(\sigma\backslash\{s,t\}; 
g(x_tu))\\
&+& \sum_{s\in\sigma\setminus t\atop s>t}(-1)^{\alpha(\sigma;s)}x_sf(\sigma\backslash\{s,t\}; 
g(x_tu))\\
&+& \sum_{s\in\sigma\setminus t\atop 
s<t}(-1)^{\alpha(\sigma;s)}\frac{x_sg(x_tu)}{g(x_sg(x_tu))}f(\sigma\backslash\{s,t\}; 
g(x_sg(x_tu)))\\
&-&\sum_{s\in\sigma\setminus t\atop 
s>t}(-1)^{\alpha(\sigma;s)}\frac{x_sg(x_tu)}{g(x_sg(x_tu))}f(\sigma\backslash\{s,t\}; 
g(x_sg(x_tu))).\\
\end{eqnarray*}
Exchanging the role of $s$ and $t$ in the second and fourth sum, 
and substituting into (\ref{b}) we obtain
\begin{eqnarray*}
(\partial\circ \psi)(f(\sigma;u))=
&-&\sum_{t\in\sigma}\sum_{s\in\sigma\setminus t\atop s<t}
      (-1)^{\alpha(\sigma;t)+\alpha(\sigma;s)}\frac{x_sx_tu}{g(x_tu)}
	  f(\sigma\backslash\{s,t\}; g(x_tu))\\
&+& \sum_{s\in\sigma}
\sum_{t\in\sigma\setminus s\atop s<t}
      (-1)^{\alpha(\sigma;t)+\alpha(\sigma;s)}\frac{x_sx_tu}{g(x_su)}
	  f(\sigma\backslash\{s,t\}; g(x_su))\\
&+& 
\sum_{t\in\sigma}\sum_{s\in\sigma\setminus t\atop s<t}
      (-1)^{\alpha(\sigma;s) +\alpha(\sigma;t)}\frac{x_{s}x_{t}u}{g(x_sg(x_tu))}
	  f(\sigma\backslash\{s,t\}; g(x_sg(x_tu)))\\
&-&
\sum_{s\in\sigma}\sum_{t\in\sigma\setminus s\atop s<t}
      (-1)^{\alpha(\sigma;s)+\alpha(\sigma;t)}\frac{x_{s}x_{t}u}{g(x_tg(x_su))}
	  f(\sigma\backslash\{s,t\}; g(x_tg(x_su))).
\end{eqnarray*}
The last two double sums in this expression cancel each other since 
by Lemma~\ref{lemma:exchange} we have  
$g(x_sg(x_tu))=g(x_tg(x_su))$ for  all $s,t\in\set(u)$. Hence a comparison with 
(\ref{a}) yields the conclusion.\end{pf}

\section{DG algebra structures on trivial extensions
			\label{section:trivialex}}
In this section we describe constructions which  in some cases allow to define algebra 
structures on free resolutions. 

We shall need the following concepts: Let $R$ be a commutative ring with a unit. 
{\em DG algebra} $A$ is a complex $(A,\partial)$ of $R$-modules with $A_i=0$ for $i<0$, 
which admits the structure of  a unitary, associative, graded commutative algebra such 
that the Leibniz rule is satisfied:
\[
\partial(ab)=\partial(a)b+(-1)^{|a|}a\partial(b)\quad \text{for all homogeneous elements}
\quad a,b\in A.
\]  
Here $|a|$ denotes the degree of $a$. 

Let $I\subset R$ be an ideal. A {\em DG algebra resolution} of $R/I$ is a 
DG algebra $A$ which is a $R$-free resolution of $R/I$. 

A {\em two-sided DG module} $M$ over $A$ is a complex of $R$-modules together with  
complex homomorphism $A\tensor_R M \to M$, $a\tensor m\mapsto am$ and $M\tensor_R A \to M$,
 $m\tensor a\mapsto ma$,  which  are unitary and associative, satisfy the Leibniz rule, 
 and the rules:
\[
(am)b=a(mb)\quad \text{and}\quad am =(-1)^{|a||m|}ma,
\]
for all homogeneous elements $a,b \in A$ and $m\in M$.  
We refer the reader to \cite{A} for details.

The following lemma which, in a different context, can be found in \cite{LA}  
and which mimics for  DG algebras the trivial extension of an algebra by a module 
(the so-called Nagata extension), is the basis of our theory.

\begin{Lemma}
\label{keylemma1}
Let $A$ be a DG algebra, $M$ a  two-sided DG A-module and 
$\psi \: M \to  A$ a DG module homomorphism.
Suppose $\psi$ satisfies the condition:
\[
\psi(m)n = m\psi(n),\quad \text{for all}\quad m, n\in M.
\]
Then the mapping cone $C(\psi)$ of $\psi$ has a DG algebra
structure with Nagata product
\begin{equation*}
     (a, m)(b, n) = (ab, (-1)^{|a|}an +mb),
\end{equation*}
for all $(a,m), (b,n) \in C(\psi)$.
\end{Lemma}
The  mapping cone  with the DG algebra structure 
as in defined in Lemma~\ref{keylemma1} will be denoted by  $A\mcone M$.
 
Let $I = (f_1,\ldots, f_{n+1}) \subset R$ be an ideal.  Set $J=(f_1,\ldots,f_n)$, 
and let $L= (f_1,\ldots, f_n):f_{n+1}$. 
Let $A$ be a free $R$-resolution of $R/J$, $M$ a free $R$-resolution of $R/L$, 
and $\psi\: M\to A$ be a complex homomorphism extending  $R/L \to  R/J$, and $C(\psi)$ 
a mapping cone.
  
Assume that $A$ and $M$ are DG algebra resolutions of $R/J$ and $R/L$, respectively. 
We want to give 
to $C(\psi)$ a DG-algebra structure 
by applying Lemma~\ref{keylemma1}. To this end, we must define a suitable action of $A$ 
on $M$,  and the complex homomorphism $\psi$
has to be chosen such that it is a DG module homomorphism over $A$ satisfying  
the condition of Lemma~\ref{keylemma1}.

We first define an action of $A$ on $M$: 
Since $J\subset L$, there is
a natural surjection
$ R/J \to  R/L$,
which induces a complex homomorphism $\phi \: A \to  M$.
Assume that  $\phi$ can be chosen to be a DG-algebra homomorphism.
Then the action of $A$ on $M$ will be defined by:
\[ 
a m = \phi(a) m\quad\text{and}\quad ma= m\phi(a) 
        \quad \text{for all}\quad a\in A\quad \text{and all} \quad m\in M,
\] 
where the product on the right hand side of the equation is multiplication
in the DG-algebra $M$. It is clear that with this action $M$ is a two-sided DG $A$-module 
over $A$. 

\begin{Lemma}
\label{lemma:prop1}
{\em (a)} Let $\psi\: M\to A$ be complex homomorphism such that 
$\phi\circ \psi =f_{n+1}\id_M$. Assume that either 
\begin{enumerate}
\item[(i)] $f_{n+1}$ is a non-zerodivisor of $R$, and  $\Im\psi$ is an ideal of $A$, or
\item[(ii)] $\phi$ is injective.
\end{enumerate}
Then $\psi$ is a DG $A$-module homomorphism.

{\em (b)} Let $\psi\: M\to A$ be a DG $A$-module homomorphism. 
The following conditions are equivalent:
\begin{enumerate}
\item[(i)] $\psi$ satisfies the condition of  Lemma~\ref{keylemma1},
\item[(ii)] $\phi\circ \psi=f_{n+1}\id_M$.
\end{enumerate}
\end{Lemma}

\begin{pf}
(a) In order to prove that $\psi$ a a DG module homomorphism, we must show that
\[
\psi(\phi(a)m) = a\psi(m), \quad \text{for all}\quad a\in A \quad \text{and}\quad m\in M.
\] 
Assuming case (i), there exists $n\in M$ with $a\psi(m)=\psi(n)$.
Also, since $\phi$ is a DG-algebra homomorphism,
we have $\phi(a\psi(m)) = \phi(a)\phi(\psi(m)) = f_{n+1}\phi(a)m$
by $\phi\circ\psi = f_{n+1}\id_{M}$.
Therefore,
\begin{eqnarray*}
f_{n+1}\phi(a)m= \phi(a\psi(m))=\phi(\psi(n))=f_{n+1}n,
\end{eqnarray*}
so that $n=\phi(a)m$. Applying $\psi$ we obtain the desired equation. 

In case (ii) the assertion follows again, since
\begin{eqnarray*}
  \phi(\psi(\phi(a)m) - a\psi(m)))
      & = & \phi(\psi(\phi(a)m)) - \phi(a\psi(m)) \\
	  & = & f_{n+1}\phi(a)m - \phi(a)\phi(\psi(m)) \\
	  & = & f_{n+1}\phi(a)m - \phi(a)(f_{n+1}m)\\
	  & = & 0.
\end{eqnarray*}

(b) Suppose  $\psi$ satisfies the condition of Lemma~\ref{keylemma1}, then 
\[
\phi(\psi(m))) = \psi(m)1_{M}=m\psi(1_{M})=mf_{n+1}=f_{n+1}m
\]
for all $m\in M$.
Hence we have $\psi\circ\phi = f_{n+1}\id_M$.
Conversely, if (ii) is satified, then for all $m, n\in M$ we have

\[
\psi(m)n =(\phi\circ\psi)(m)n=f_{n+1}mn= m(f_{n+1}n)=
m(\phi\circ\psi)(n) = m\psi(n).
\]
\end{pf}

\section{Koszul type resolutions
  				\label{section:koszultype}}
We introduce the following notion. 

\begin{Definition}
\label{def}
{\em Let $I\subset R$ be an ideal. 
A DG algebra resolution $A$ of 
$R/I$ over $R$ is of {\em Koszul type} (of length $n$), if for all $i$
\begin{enumerate}
\item [(i)] $\rank A_i={n\choose i}$,
\item [(ii)] the homomorphism $A_i\to \Hom_R(A_{n-i},A_n)$ which assigns 
       to each $a\in A_i$ the map $\alpha_a\: A_{n-i}\to A_n$ with $\alpha_a(b)=ab$ 
	   is injective. 
\end{enumerate}}
\end{Definition}

Let, as in Section~\ref{section:trivialex}, $A$ and $M$ be DG algebra resolutions of 
$R/J$ and $R/I$, respectively, and let $\phi\:A\to M$ be a DG algebra homomorphism.

\begin{Lemma}
\label{lemma:ktype}
Assume that $R$ is a domain with quotient field $Q$, and that $A$ and $M$ are 
of Koszul type of length $n$. Fix a basis $e$ for $A_n$ and a basis $\bar{e}$ for $M_n$, 
and assume that $\phi(e)=\delta\bar{e}$ with $\delta\neq 0$. 
Then there is a unique complex homomorphism
$\tilde{\phi}\: M\tensor_RQ\to A\tensor_RQ$ satisfying
\[
\tilde{\phi}(m)a=m\phi(a),
\]
for all $0\leq i\leq n$, $a\in A_{n-i}$ and $m\in M_i$, 
where the equation means equality between the coefficients of the bases $e$ and $\bar{e}$.
Moreover, $\tilde{\phi}\circ \phi=\delta \id_{A\tensor Q}$ 
and $\phi\circ \tilde{\phi}=\delta\id_{M\tensor Q}$. 
\end{Lemma}

\begin{pf} Let $\epsilon\: M_n\tensor Q\to A_n\tensor Q$ be the isomorphism with 
$\epsilon(\bar{e})=e$. Given $m\in M_i\tensor Q$, we define the $Q$-linear map 
$\gamma\: A_{n-i}\tensor Q\to A_n\tensor Q$ by $\gamma(a)=\epsilon(m\phi(a))$. 
(Here we write $\phi$ instead of $\phi\tensor Q$, in order to simplify notation.) 

Since by assumption the natural map $A_i\to \Hom_R(A_{n-i},A_n)$ is injective, 
the induced map 
$A_i\tensor Q\to \Hom_Q(A_{n-i}\tensor Q,A_n\tensor Q)$ 
is again injective, and hence must even be bijective since it is a linear map of 
vector spaces of equal dimension. 
Therefore, there exists $b\in A_i\tensor Q$ with $ba=\gamma(a)=m\phi(a)$ 
for all $a \in A_{n-i}$, and we set $\tilde{\phi}(m)=b$. 
Then $\tilde{\phi}(m)a=m\phi(a)$.

For arbitrary elements
$m\in M_{i+1}\tensor Q$ and $a\in A_{n-i}\tensor Q$, we have 
\begin{equation*}
\tilde{\phi}(\partial(m)) a  =  \partial(m)\phi(a).
  \end{equation*}
Since $m\phi(a) \in M_{n+1}\tensor Q = 0$,
we have $0 = \partial(m\phi(a)) = \partial(m)\phi(a) + (-1)^{\mid m\mid}m\partial(\phi(a))$.
Hence 
\begin{eqnarray*}
\partial(m)\phi(a) &=& -(-1)^{\mid m\mid}m\partial(\phi(a)) \\
                   &=& -(-1)^{\mid m\mid}m\phi(\partial(a)) \\
				   &=& -(-1)^{\mid m\mid}\tilde{\phi}(m)\partial(a),
\end{eqnarray*}
and since $\tilde{\phi}(m) a \in A_{n+1} = 0$, 
we have $0 = \partial(\tilde{\phi}(m)a)
           = \partial(\tilde{\phi}(m))a 
		   	 + (-1)^{\mid \tilde{\phi}(m)\mid}\tilde{\phi}(m)\partial(a)$,
so that 
\begin{equation*}
\tilde{\phi}(m)\partial(a) = (-1)^{\mid m\mid + 1}\partial(\tilde{\phi}(m))a.
\end{equation*}
Consequently, 
\begin{equation*}
     \tilde{\phi}(\partial(m))a = \partial(\tilde{\phi}(m))a 
\end{equation*}
for all $a\in A_{n-i}\tensor Q$ and all $m\in M_{i+1}\tensor Q$. Therefore, 
since $A$ is of Koszul type we conclude that 
\begin{equation*}
     \tilde{\phi}(\partial(m)) = \partial(\tilde{\phi}(m)),
\end{equation*}
which shows that
$\tilde\phi$ is indeed a complex homomorphism.

Let $a\in A_i\tensor Q$ and $b\in A_{n-i}\tensor Q$. Then  $\tilde{\phi}(\phi(a))b =
\phi(a)\phi(b) = \phi(ab) = \delta ab$.
Hence since $A$ is of Koszul type, it follows that $\tilde{\phi}(\phi(a)) = \delta a$. 
In other words, $\tilde{\phi}\circ\phi = \delta \id_{A\tensor Q}$.
This implies in particular that $\phi\: A_i\tensor Q\to M_i\tensor Q$ is injective. 
However, since $A_i\tensor Q$ and $M_i\tensor Q$ have the same $Q$-dimension, 
we see that $\phi$ is an isomorphism with $\phi^{-1}=\delta^{-1}\tilde{\phi}$. 
Therefore, $\phi\circ \delta^{-1}\tilde{\phi}=\id_{M_i\tensor Q}$, 
and so $\phi\circ \tilde{\phi}=\delta\id_{M_i\tensor Q}$, as desired. 
\end{pf}

Now Lemma~\ref{lemma:prop1} and Lemma~\ref{lemma:ktype} imply immediately

\begin{Corollary}
\label{cor:psimap}
With  the assumptions and the notation of Lemma~\ref{lemma:ktype},
suppose that $f_{n+1}\tilde{\phi}(M)\subset\delta A$, and set $\psi=
(f_{n+1}/\delta)\tilde{\phi}$. If $M$ is viewed a two-sided $A$-module via
$\phi$, then $\psi\:M\to A$ is a DG module homomorphism satisfying the condition 
of Lemma~\ref{keylemma1}. In particular, $A*M$ is defined, 
and $A*M$ is of Koszul type of length $n+1$.
\end{Corollary}

\begin{pf}
It remains to be shown that $A*M$ is of Koszul type of length $n+1$. In fact, we have 
$(A*M)_i=A_i\dirsum M_{i-1}$, and so 
$\rank (A*M)_i=\rank A_i+\rank M_{i-1}={n\choose i}+{n\choose i-1}={n+1\choose i}$. 

Finally, let $(a,m)\in (A*M)_i$, $(a,m)\neq 0$. We must show that there exists 
$(b,n)\in (A*M)_{n+1-i}$ such that $(a,m)(b,n)\neq 0$. 
By the definition of our multiplication we get
$(a,m)(b,n)=(0,(-1)^{|a|}\phi(a)n+m\phi(b))$. 

There are two cases to consider. If $a\neq 0$, we let $b=0$, and have to show 
that there exists $n\in M$ with $\phi(a)n\neq 0$. 
Now since $M$ is of Koszul type, there exists $w\in M_i$ such that $wn\neq 0$. 
But recall that $A_i\tensor Q\to M_i\tensor Q$ is an isomorphism via $\phi$. 
Thus, since $a\neq 0$, $\phi(a)n\neq 0$ for arbitrary $n\neq 0$, and 
by the isomorphism $A_n\tensor Q\to M_n\tensor Q$ between one dimensional
vector spaces we have $\lambda\in Q$, $\lambda\neq 0$ such that
$\phi(a)n = \lambda wn \neq 0$.

In the second case, $m\neq 0$, and we let $n=0$. Then we have to find $b\in A_{n+1-i}$ 
such that $m\phi(b)\neq 0$. The rest of the argument is the same as in 
the first case.
\end{pf}

\begin{Definition}
\label{def2}
{\em A sequence $f_1, \ldots, f_n$ in $R$ is called a {\em Koszul sequence}, 
if for all $i=1,\ldots, n$
\begin{enumerate}
\item[(i)] $R/(f_1,\ldots,f_i)$ has a Koszul type resolution  $A^{(i)}$ of length $i$;
\item[(ii)] $R/((f_1,\ldots,f_{i-1}):f_{i})$  has a Koszul type resolution  $M^{(i-1)}$ 
of of length $i$; 
\item[(iii)] $A^{(i)}\iso A^{(i-1)}*M^{(i-1)}$.
\end{enumerate}}
\end{Definition}

We will now consider some examples.

\begin{Example}
\label{example:reg}
Regular sequences are Koszul sequences. 
\end{Example}

\begin{pf}
Let $f_1,\ldots,f_n$ be regular sequence. For a given $i$ we let $A^{(i)}$ be t
he Koszul complex for the sequence $f_1,\ldots,f_i$. 
Since $(f_1,\ldots,f_i):f_{i+1}=(f_1,\ldots,f_i)$, we may choose $M^{(i)}=A^{(i)}$, 
and  $\phi=\id_{A{(i)}}$. Then $\psi=f_{i+1}\id_{A^{(i)}}$. It is then easy to see that 
$A^{(i+1)}=A^{(i)}* A^{(i)}$ is the Koszul complex for the sequence $f_1,\ldots,f_{i+1}$.
\end{pf}

\begin{Example}
\label{example:taylor}
Monomial sequences are Koszul sequences.
\end{Example}

\begin{pf}
Let $f_1,\ldots,f_n$ be a sequence of monomials. 
The Taylor resolution $T=T(f_1,\ldots,f_n)$ of this sequence 
admits a natural DG algebra structure, as shown by Gemeda \cite{G} (see also \cite{F}). 
As an $R$-module $T_k$ is the $k$th exterior power of $T_1= \Dirsum_{i=1}^n Re_i$ with 
$R$-basis $\{e_\sigma\: \sigma\subset [n], |\sigma|=k\}$, 
where $e_\sigma=e_{i_1}\wedge e_{i_2}\wedge \cdots\wedge e_{i_k}$ 
for  $\sigma=\{i_1<i_2<\cdots<i_k\}$.  
The chain map $\partial\: T_i\to T_{i-1}$ is defined by 
\[
\partial(e_\sigma)=\sum_{i\in \sigma}(-1)^{\sigma(\sigma,i)}
							  \frac{f_\sigma}{f_{\sigma\setminus \{i\}}}
							  e_{\sigma\setminus \{i\}},
\]
where for $\tau\subset [n]$, $f_\tau$ denotes the least common multiple of 
the monomials $f_i$ with $i\in \tau$, and where $\sigma(\sigma,i)=|\{j\in \sigma\: j<i\}|$. 
According to Gemeda, the DG algebra structure on $T$ is given by
\[
e_\sigma e_\tau=\frac{f_\sigma f_\tau}{f_{\sigma\union \tau}}e_\sigma\wedge e_\tau.
\]

It is also known, and easy to see, that $T(f_1,\ldots,f_n)$ is obtained as the mapping cone
of $\psi\:T(g_1,\ldots,g_{n-1}) \to T(f_1,\ldots, f_{n-1})$, where $g_i=f_i/[f_i,f_n]$ for 
$i=1,\ldots,n$, and where $\psi(\bar{e}_\sigma)=(f_{\sigma\union \{n\}}/f_\sigma)e_\sigma$.
 Here $\{\bar{e}_\sigma\: \sigma\subset [n-1]\}$ denotes the natural basis of 
 $T(g_1,\ldots,g_{n-1})$. 

We now define an $R$-module homomorphism $\phi\: T(f_1,\ldots,f_n)\to T(g_1,\ldots,g_n)$
by
\[
\phi(e_\sigma)=\frac{f_\sigma f_{n+1}}{f_{\sigma\union\{n+1\}}}\bar{e}_\sigma,
\quad\text{for all}\quad \sigma\subset [n].
\]
It is easy to check that $\phi$ is an injective DG algebra homomorphism, and that 
$\phi\circ \psi=f_{n+1}\id$. Therefore, 
Lemma~\ref{lemma:prop1} implies that $\psi$ is a DG $A$-module homomorphism satisfying the 
condition
Lemma~\ref{keylemma1}. Thus the DG algebra  $T(f_1,\ldots,f_n)*T(g_1,\ldots,g_n)$ 
is defined.

Consider the $R$-module homomorphism
\[
\alpha: T(f_1,\ldots, f_{n+1})\To T(f_1,\ldots,f_n)*T(g_1,\ldots,g_n)
\]
with
\[
\alpha(e_\sigma)=\left\{\begin{array}{ll}
(e_\sigma,0),&\text{if $n+1\not\in \sigma$},\\
(0,\bar{e}_{\sigma\setminus \{n+1\}}), &\text{if $n+1\in \sigma$}.
\end{array}
\right.
\]
We leave it to the reader to check that this DG algebra isomorphism.

It may be worthwhile to notice that with the notation of Corollary~\ref{cor:psimap} one has
$\psi= (f_{n+1}/\delta)\tilde{\phi}$ this case $\delta=(f_{[n]}f_{n+1}/f_{[n+1]})$.

\end{pf}

\begin{Example}
\label{example:almost}
Let $f_1,\ldots, f_n$ and $g_1,\ldots, g_n$ be regular sequences, such that 
\[
(f_1,\ldots,f_n)\subset (g_1,\ldots,g_n).
\] 
Let $f_i=\sum_{j=1}^na_{ij}g_j$ for $i=1,\ldots,n$, and set $\Delta=\det(a_{ij})$. 
Then the almost complete intersection $(f_1,\ldots,f_n,\Delta)$ is a Koszul sequence.
\end{Example}

\begin{pf}
The initial sequence $f_1,\ldots,f_n$ is a Koszul sequence by 
Example~\ref{example:reg}.
 Next we observe observe that $(g_1,\ldots,g_n)=(f_1,\ldots,f_n):\Delta$. We let $A^{(n)}$ 
 be the Koszul complex of the sequence $f_1,\ldots, f_n$,  and $M^{(n)}$ the Koszul complex
  of the sequence $g_1,\ldots,g_n$. Let $e_1,\ldots, e_n$ be the $R$-module basis of 
  $A^{(n)}_1$  with $\partial(e_i)=f_i$ for $i=1,\ldots,n$, and $h_1,\ldots,h_n$ the 
  $R$-module basis of  $M^{(n)}_1$  with $\partial(h_i)=g_i$ for $i=1,\ldots,n$. 
  Then the unique algebra homomorphism $\phi\: A^{(n)}\to M^{(n)}$ with 
  $\phi(e_i)=\sum_{j=1}^na_{ij}h_j$ for $i=1,\ldots,n$ extends the epimorphism 
  $R/(f_1,\ldots,f_n)\to R/(g_1,\ldots,g_n)$, and
\[
\phi(e_1\wedge\cdots\wedge e_n)=\Delta (h_1\wedge\cdots\wedge h_n).
\]
Thus Corollary~\ref{cor:psimap} implies that the mapping cone $C(\psi)=A^{(n)}*M^{(n)}$  
with  $\psi=\tilde{\phi}$ is of Koszul type.
\end{pf}   

Note that the almost complete intersection considered in Example~\ref{example:almost} 
is directly linked to complete intersection. More generally, let $I\subset R$ be 
a perfect ideal of grade $g$ in a Gorenstein ring $R$, $L\subset I$ a complete intersection 
ideal of the same grade, and $J=L:I$ the linked ideal. Then the canonical module $\omega_A$
 of $A$ is isomorphic to $\Hom_R(R/I,R/L)$ (see for example \cite{BH}), 
 and $\Hom_R(R/I,R/L)\iso (L:I)/L=J/L$. Therefore one obtains an exact sequence 
\[
0\To \omega_A\To R/L\To R/J\To 0.
\]
The $R$-dual $A^*=\Hom_R(A,R)$ (with $A^*_i=\Hom_R(A_{g-i},R)$) is a graded minimal free 
$R$-resolution of $\omega_A$ (cf.\ \cite{BH}), and since $K$  is self dual, we may lift 
the $R$-module homomorphism $\omega_A\to R/L$ to a graded complex homomorphism 
$\psi\: A^*\to K^*$. Then the  mapping cone $C(\psi)$ is a graded free resolution of $R/J$,
 as is well-known.

In case $R/I$ is Gorenstein, in which case $A\iso A^*$, one could hope to define an 
algebra structure on $C(\psi)$ just as in
Example~\ref{example:almost}. But this is not possible since we would need that the 
composition $A\to K\to A$ is the multiplication map by an element of $R$. By rank reasons 
this could only be possible if $\rank A_i=\rank K_i= {d\choose i}$ for all $i$.

On the other hand, if we suppose  that $A$ is two-sided DG $K$-module and that the 
epimorphism $R/L\to R/I$ can be extended to DG $K$-module homomorphism $\phi\: K\to A$, 
then $\psi\: A^*\to K^*$ can be chosen such that $C(\psi)$ has a natural two-sided 
DG $K$-module structure. 

In fact, we first define the structure of a two-sided  DG $K$-module  on $A^*$ as follows:
 For $\alpha\in A^*_{i}$ and $c\in K_j$ we let $c\alpha, \alpha c\in A^*_j$  with 
 $c\alpha(a)=\alpha(a\phi(c))$ and $\alpha c(a)=\alpha(\phi(c)a)$ for all $a\in A_{n-i-j}$. 
Then let $\psi=\phi^*$, the $R$-dual of $\phi$. It is then easily checked that 
$\psi\: A^*\to K^*$ is a DG $K$-module homomorphism. 
This immediately implies that $C(\psi)$ is a two-sided DG $K$-module.

\end{document}